\documentclass[a4paper 12pt]{article}
\usepackage{latexsym,amsmath, amsthm, amssymb}
\newtheorem{lem}{Lemma}[section]
\newtheorem{thm}{Theorem}[section]
\newtheorem{definition}{Definition}[section]

\newtheorem{remark}{Remark}[section]
\newtheorem{example}{Example}[section]
\newtheorem{pro}{Proposition}[section]

\def\bc{\begin{center}}
\def\ec{\end{center}}

\begin{document}
 \abovedisplayskip=8pt plus 1pt minus 1pt
\belowdisplayskip=8pt plus 1pt minus 1pt
\thispagestyle{empty} \vspace*{-3.0truecm} \noindent
\parbox[t]{6truecm}{\footnotesize\baselineskip=11pt\noindent  {} 
 } \hfill
\vspace{1 true cm}

\bc{\large\bf Continuous dependence property of BSDE with constraints\footnote{The work of Helin Wu is  Supported by the Scientific and Technological Research
Program of Chongqing Municipal Education Commission  (KJ1400922). The work of Yong Ren is supported by the
National Natural Science Foundation of China (11371029). The work of Feng Hu is supported by the
National Natural Science Foundation of China (11301295), the Program Foundation of Ministry of Education of
China (20123705120005) and the  Program  for
Scientific  Research  Innovation Team in
Applied Probability and  Statistics
of  Qufu Normal University (0230518).}}
 \ec

\bc{\bf Helin Wu  \\
{\small\it School of Mathematics,  Chongqing  University of Technology, Chongqing 400054, China\\
\small\it \quad Corresponding author. e-mail: wuhewlin@gmail.com}}\ec

\bc{\bf Yong Ren  \\
{\small\it Department of Mathematics, Anhui Normal University, Wuhu 241000, China\\
\small\it \quad  e-mail: renyong@126.com}}\ec

\bc{\bf Feng Hu  \\
{\small\it School of Statistics, Qufu Normal University, Qufu 273165, China\\
\small\it \quad  e-mail: hufengqf@163.com}}\ec

\begin{abstract}
   In this paper, we study continuous properties of  adapted solutions for  backward stochastic differential equations
   with constraints (CBSDEs in short). Comparing with many
    existing literatures about this topic,
   our case is very general in the sense that constraints are formulated by  general non-negative real functions.
   In  general  case, we proved a continuous property from below and a
lower semi-continuous property of the minimal super-solution of
CBSDE in
   its effective domain. Furthermore,
   in the special convex case, we obtained a continuous property with the help of convex analysis.

\end{abstract}
{\bf Keywords:} backward stochastic differential
equation with constraint (CBSDE), continuous dependence, convex functional,
minimal solution.

\section {Introduction}
By Pardoux and Peng \cite{pepsg}, we know that there exists a unique
adapted and square integrable solution to a backward stochastic differential equation (BSDE in short) of the type
$$y_t=\xi+\int_t^T g(s,y_s,z_s)ds-\int_t^Tz_sdW_s,\quad
0\leq t\leq T,\eqno(1.1)$$ under some suitable conditions on the
function $g$ and the terminal value $\xi$.

It is useful for this kind of BSDE to calculate pricing and hedge
claims in mathematical finance,  where $y(t)$ represents the
wealth process and $z(t)$ represents the portfolio process, see
Karoui et al. \cite{kpq}. However,  the information of financial
markets may be incomplete  or
 suffers other constraints. To calculate pricing and hedge claims with constraints on wealth and portfolio,
  BSDE with constraint (CBSDE in short) was introduced, that is
$$y_t=\xi+\int_t^Tg(s,y_s,z_s)ds-\int_t^Tz^*_sdW_s+\int_t^TdC_s, \quad 0\leq t \leq T, \eqno(1.2)$$
$$(y(t),z(t))\in \Gamma_t,\quad a.e., \quad \text{on} \quad [0,T]\times \Omega, \eqno(1.3)$$
where $C(t)$ is an   RCLL increasing process
 with $C(0)=0$, $\Gamma_t:=\{(y,z)|\phi(t,y,z)=0\}\subset \mathbb{R}\times\mathbb{R}^d$ and
$\phi(t,y,z): [0,T]\times \mathbb{R}\times \mathbb{R}^d\rightarrow
\mathbb{R}^+$.

Among many properties of the solutions for BSDEs, the continuous
dependence property is very important in the theory of BSDEs. In
the classic case with no constraints, it is easy to derive the
continuous dependence property of $y(t)$ with the terminal value
$\xi$ under the conditions that the  generator $g$ is Lipschitz in
both variables $y$ and $z$, and $(g(t, 0, 0))_{t\in[0,T ]}$ are
square integrable. One can see Pardoux and Peng \cite{pepsg}.
 When the generator $g$ of BSDEs does not satisfy the above conditions, the continuous dependence property is weakened
generally. In this respect, we can see Fan et al. \cite{fwz} for
more details.  Recently, Drapeau et al. \cite{dhk} discussed the
continuous property of the minimal solutions of  BSDEs with convex
generators satisfying some assumptions in a more general space
$L_T^0(\mathbb{R})$, the set of $\mathcal{F}_T$-measurable
functions. The origin to consider the minimal solutions of CBSDEs
comes from the super-hedging and super-pricing problem in
incomplete market or constrained market, one can see Karatzas and
Shreve \cite{kiss} and the references therein. To our best
knowledge, although many works have been done  on CBSDEs, the
continuous dependence property was ignored  for a long time. In
this paper, our aim is to study the continuous properties of
CBSDEs.
 We firstly prove some  continuous properties in semi sense about
$g_\Gamma$-solutions  with general coefficients. Then, when both $g$
and $\phi$ are convex, we prove a continuous property in the
interior of its  domain.

 This paper is organized as follows. In Section 2, we state the
 framework in  Peng \cite{psg} and collect some propositions about
 $g_\Gamma$-solution. The main results are obtained in Section 3 and we give a counter-example to
 show that $g_\Gamma$-solution is not continuous in non-convex case. Some key results about convex analysis
 used in our paper are listed in Section 4.

\section {BSDE  and $g_\Gamma$-solution of CBSDE}

In what follows, let $T>0$ be a fixed deterministic terminal time
and $\{W_t\}_{t\geq0}$ be a $d$-dimensional standard Brownian
motion defined on a complete probability space $(\Omega, \mathcal
{F}, P)$. We denote the natural filtration generated by
$\{W_t\}_{0\leq t\leq T}$ and augmented by all $P$-null sets by
$\mathbb{F}:=\{\mathcal {F}_s, 0\leq s\leq T\}$. Let
$L_T^2(\mathbb{R})$ denote the space of ${\cal F}_T$-measurable
random variables $\xi$ satisfying $E|\xi|^2<\infty$. Suppose a
function $g(\omega,t,y,z):\Omega\times[0,T]\times \mathbb{R}
\times \mathbb{R}^d\rightarrow\mathbb{R}$ satisfies the  uniformly
Lipschitz condition, i.e., there exists a positive constant $M$
such that for all $(y_1,z_1),(y_2,z_2)\in \mathbb{R}
\times\mathbb{R}^d$
$$|g(\omega,t,y_1,z_1)-g(\omega,t,y_2,z_2)|\leq M(|y_1-y_2|+|z_1-z_2|)
\eqno(A1)$$ and
$$g(\cdot,0,0)\in H_T^2(\mathbb{R}), \eqno(A2)$$
where $H_T^2(\mathbb{R})$ denotes  the space of predictable
processes $\varphi:\Omega\times[0,T]\rightarrow \mathbb{R}$
satisfying $\parallel
\varphi\parallel^2=E\int_0^T|\varphi(s)|^2ds<+\infty.$ By Pardoux
and Peng \cite{pepsg}, we know that if $g$ satisfies $(A1)$ and
$(A2)$, for any   $\xi\in L_T^2(\mathbb{R})$, the BSDE (1.1) has a
unique adapted and square integrable solution
$(y(t),z(t))_{t\in[0,T]}$. We call $y(t)$ as $g$-solution.

Sometimes we also need  the following assumption:
$$g(\cdot,\cdot,0)=0.\eqno(A3)$$

When constraint conditions as in (1.3) is considered,
$g$-super-solutions are defined as (1.2) for $\xi\in
L_T^2(\mathbb{R})$ and  $g$ satisfying conditions $(A1)$ and
$(A2)$. Among all such $g$-super-solutions, we mainly concern
about the minimal one defined as below.

\begin{definition}\rm ($g_\Gamma$-solution) A $g$-super-solution $(y_t, z_t, C_t)$ is said to be
 the minimal solution, given $y_T=\xi$,
subjected to the constraint $(1.3)$ if for any other
g-super-solution $(y'_t, z'_t, C'_t)$ satisfying $(1.3)$  with
$y'_T=\xi$, we have $y_t\leq y'_t $ a.e.. We call the  minimal
solution  as $g_\Gamma$-solution and denote it by
$\mathcal{E}_t^{g,\phi}(\xi)$
 for convenience when $g,\phi$ satisfy $(A3)$ for any $t\in [0,T]$.
\end{definition}
For any $\xi\in L^2_T(\mathbb{R})$, we denote
$\mathcal{H}^{\phi}(\xi)$  the set of g-super-solutions
$(y_t,z_t,C_t)$ subjecting to $(1.3)$ with $y_T=\xi$.  When
$\mathcal{H}^{\phi}(\xi)$ is not empty, the existence of the
minimal solution, namely  $g_\Gamma$-solution was proved in
 Peng \cite{psg} by some kind of penalization method.

The convexity of $\mathcal{E}_t^{g,\phi}(\xi)$ can be easily deduced
from the convexities of $g$ and $\phi$.

\begin{pro}\label{con1}
Suppose that $\phi(\omega,t,y,z)$ : $\Omega\times[0,T]\times
\mathbb{R}\times \mathbb{R}^d\rightarrow \mathbb{R}^+$ and
$g(\omega,t,y,z)$ : $\Omega\times[0,T]\times\mathbb{R}\times
\mathbb{R}^d\rightarrow\mathbb{R}$ are convex functions satisfying
the assumptions $(A1)$ and $(A3)$, then
$$\mathcal{E}_t^{g,\phi}(a\xi+(1-a)\eta)\leq a\mathcal{E}_t^{g,\phi}(\xi)+(1-a)\mathcal{E}_t^{g,\phi}(\eta)
\quad a.e., \quad\text{on} \quad[0,T]\times\Omega$$
 for any $\xi,\eta\in L_T^2(\mathbb{R})$  and $a\in[0,1]$ when $\mathcal{H}^{\phi}(\xi), \mathcal{H}^{\phi}(\eta)$ are not empty.
\end{pro}
\proof
When $\mathcal{H}^{\phi}(\xi), \mathcal{H}^{\phi}(\eta)$ are not empty,
it is easy to see that  $\mathcal{H}^{\phi}(a\xi+(1-a)\eta)$ is not empty for $a\in[0,1]$.

According to  Peng \cite{psg}, the solutions
$\{y_t^m(\xi),m=1,2,\cdots\}$ of
$$y_t^m(\xi)=\xi + \int_t^Tg(s,y_s^m(\xi),z_s^m)ds+A_T^m-A_t^m-\int_t^T(z_s^m)^*dW_s.$$
is an increasing sequence  converging  to
$\mathcal{E}_t^{g,\phi}(\xi)$, where $$A_t^m: =
m\int_0^t\phi(s,y_s^m,z_s^m)ds.$$ For any fixed $m$, by the
convexity of $g$ and $\phi$, $y_t^m(\xi)$ is  convex in $\xi$,
that is
$$y_t^m(a\xi+(1-a)\eta)\leq ay_t^m(\xi)+(1-a)y_t^m(\eta),$$
taking limit on both sides as $m\rightarrow\infty$, we get the
required result.
\hspace*{\fill}$\Box$\\

By the penalization method (see Peng \cite{psg}), we can easily
obtain the following comparison theorem of
$\mathcal{E}_t^{g,\phi}(\xi)$ .
\begin{pro}\label{com}
Under the same assumptions as Proposition 2.1, we have
$$\mathcal{E}_t^{g,\phi}(\xi )\leq
\mathcal{E}_t^{g,\phi}(\eta)$$ for any $\xi, \eta\in L_T^2(\mathbb{R})$ with
$P(\eta\geq\xi)=1$ when $\mathcal{H}^{\phi}(\xi), \mathcal{H}^{\phi}(\eta)$ are not empty.
\end{pro}

\section {Continuous property about $g_\Gamma$-solution }
In this section, we firstly prove a continuous property from below and a
lower semi-continuous property of  $\mathcal{E}_t^{g,\phi}(\xi)$.
 Then, we prove a continuous property in the interior of the domain
 when the coefficients are convex.

\begin{thm}\label{cdp}
Suppose  $\phi(\omega,t,y,z)$ : $\Omega\times[0,T]\times \mathbb{R}\times \mathbb{R}^d\rightarrow \mathbb{R}^+$
and $g(\omega,t,y,z)$ : $\Omega\times[0,T]\times \mathbb{R}\times \mathbb{R}^d\rightarrow \mathbb{R}$
 satisfy the assumptions $(A1)$ and $(A3)$, $\xi_n\leq \xi, \, n=1,2,\cdots $ a.s. and $E|\xi_n-\xi|^2\rightarrow 0$ as  $n\rightarrow \infty$, then
  $$E|\mathcal{E}_t^{g,\phi}(\xi_n)-\mathcal{E}_t^{g,\phi}(\xi)|^2\rightarrow
 0$$
if   $\mathcal{H}^{\phi}(\xi)$ is not empty.
 \end{thm}
 \proof Firstly, it is easy to see that $\mathcal{H}^{\phi}(\xi_n)$ is not empty due to the assumptions of $\xi_n\leq \xi, \forall n$ a.s. and the nonempty of $\mathcal{H}^{\phi}(\xi)$. In fact, if $(y_t,z_t,C_t)$ is a super-solution with terminal value $\xi$, then
$(y'_t,z_t,C'_t)$ is a super-solution with terminal value $\xi_n$, where
\[
C'_t=\left\{
\begin{array}{ll}
C_t &\text{when}\quad 0\leq t<T;\\
C_T+\xi-\xi_n &\text{when}\quad t=T.
\end{array}
\right.
\]

\[
y'_t=\left\{
\begin{array}{ll}
y_t &\text{when}\quad 0\leq t<T;\\
\xi_n &\text{when}\quad t=T.
\end{array}
\right.
\]

Let us  consider the functional $\varphi(\eta):=
E|\mathcal{E}_t^{g,\phi}(\eta)-\mathcal{E}_t^{g,\phi}(\xi)|^2$ on
the convex set $\tilde{K}=\{\eta\in D|\eta\leq \xi \quad a.s.\}$ for a
fixed $\xi\in L_T^2(R)$.

Taking any sequence $\{\eta_n\in\tilde{K},n=1,2,\cdots\}$  converges
to $\eta\in\tilde{K}$ a.s., we want to   prove  that
$\varphi(\eta)\geq c$ if $\varphi(\eta_n)\geq c$  for  all  $n$. To this aim,
 let $y_t^m(\eta_n)$ be  the  approximating  sequence  of
$\mathcal{E}_t^{g,\phi}(\eta_n)$  as  in  Proposition 2.1, i.e.,
$\{y_t^m(\eta_n),m=1,2,\cdots\}$ is a sequence dominated by
$\mathcal{E}_t^{g,\phi}(\eta_n)$ and converges increasingly to
$\mathcal{E}_t^{g,\phi}(\eta_n)$  as $m\rightarrow\infty$.  Since   $\eta_n\leq\xi$ a.s. and $y_t^m(\eta_n)\leq
\mathcal{E}_t^{g,\phi}(\eta_n)\leq\mathcal{E}_t^{g,\phi}(\xi)$ for
any $n$ and $m$,   it
is obvious that
$E|y_t^m(\eta_n)-\mathcal{E}_t^{g,\phi}(\xi)|^2\geq
E|\mathcal{E}_t^{g,\phi}(\eta_n)-\mathcal{E}_t^{g,\phi}(\xi)|^2\geq
c$. For any fixed $m$, by the continuous dependence
property of unconstrained BSDE,
$E|y_t^m(\eta_n)-y_t^m(\eta)|^2\rightarrow 0$ as
$n\rightarrow\infty$, which shows that
$E|y_t^m(\eta)-\mathcal{E}_t^{g,\phi}(\xi)|^2\geq c$  for any  $m$. Furthermore,
$E|\mathcal{E}_t^{g,\phi}(\eta)-\mathcal{E}_t^{g,\phi}(\xi)|^2\geq
c$ by  monotone convergence theorem.

Set $\beta_n:=\varphi(\xi_n)=E|\mathcal{E}_t^{g,\phi}(\xi_n)-\mathcal{E}_t^{g,\phi}(\xi)|^2$. Suppose  on the contrary that $\beta_n\nrightarrow 0$ as
$n\rightarrow \infty$,  then there exists some subsequence of
$\{\xi_n,n=1,2,\cdots\}$ (for convenience, we still
denote it as $\{\xi_n,n=1,2,\cdots\}$)
 such that $\beta_n=
E|\mathcal{E}_t^{g,\phi}(\xi_n)-\mathcal{E}_t^{g,\phi}(\xi)|^2\geq\delta$
for some $\delta>0$.  If we take  $\eta_n$ as $\xi_n$ converges to $\eta=\xi \quad a.s.$ in
the above argument, then there will  be a  contradiction  since $\varphi(\eta):=
E|\mathcal{E}_t^{g,\phi}(\eta)-\mathcal{E}_t^{g,\phi}(\xi)|^2=
E|\mathcal{E}_t^{g,\phi}(\xi)-\mathcal{E}_t^{g,\phi}(\xi)|^2=0\geq \delta$.

\hspace*{\fill}$\Box$

\begin{remark}\rm Let $L_T^\infty(\mathbb{R})$ denote the
space of all P-essentially bounded random variables. We define a risk measure
 $\rho(\cdot):=\mathcal{E}_0^{g,\phi}(-\cdot)$ on $L_T^\infty(\mathbb{R})$. It is interesting to note that according to  Peng and Xu \cite{psgxmy1},  and Rosazza Gianin \cite{ger}, when $g(t,y,z)$ and  $\phi(t,y,z)$
are independent of $y$ and convex in $z$, Theorem 3.1 implies  that $\rho(\cdot)$ satisfies the important Fatou property (the Fatou property can be seen in F\"{o}llmer and Schied  \cite{fhsa}).

\end{remark}
\begin{remark}\rm
In Drapeau et al. \cite{dhk}, a continuous dependence property has been proved when $\{\xi_n\in L_T^2(\mathbb{R}), \ n=1,2,\cdots\}$
  converging  to $\xi\in L_T^2(\mathbb{R})$ increasingly. This is a monotonic continuous  property. But in Theorem \ref{cdp}, we only assume
  $E|\xi_n-\xi|^2\rightarrow 0$ as $n\rightarrow \infty$ and $\xi_n\leq \xi \quad a.s.$ . We call this property the continuous dependence property from below.
\end{remark}

In the following, we  prove the lower semi-continuity of  $\mathcal{E}_0^{g,\phi}(\cdot)$.
Let $D:=\{\xi\in
L_T^2(\mathbb{R})|-\infty<\mathcal{E}_0^{g,\phi}(\xi)<\infty\}$ denote the effective
 domain of $\mathcal{E}_0^{g,\phi}(\cdot)$. Up to now,
the structure of $D$ is not clear.
 In this paper, we assume that
 $$ D \quad \text{is norm closed in} \quad L_T^2(\mathbb{R}). \eqno(A4)$$

Then  we have the following result.
\begin{thm}\label{lcdp}
Suppose  $\phi(\omega,t,y,z)$ : $\Omega\times[0,T]\times \mathbb{R}\times \mathbb{R}^d\rightarrow \mathbb{R}^+$
and $g(\omega,t,y,z)$ : $\Omega\times[0,T]\times \mathbb{R}\times \mathbb{R}^d\rightarrow \mathbb{R}$
satisfy the assumptions $(A1)$ and $(A3)$. If $(A4)$ holds, then for any $k\in \mathbb{R}$, $A_k$ is
closed  in $L_T^2(\mathbb{R})$, where $A_k:=\{\xi\in D|\mathcal{E}_0^{g,\phi}(\xi)\leq k\}$, the $k$-level set of $\mathcal{E}_0^{g,\phi}(\xi)$.
\end{thm}
\proof Suppose a sequence $\{\xi_n,n=1,2\cdots\}\subset A_k$
converges to some $\xi\in L_T^2(\mathbb{R})$. By the closeness of
$D$, $\xi\in D$. For any $\xi_n$, we take
$\{y_0^m(\xi_n),m=1,2,\cdots\}$ as in Proposition \ref{con1}. Since
$y_0^m(\xi_n)$ converges increasingly to
$\mathcal{E}_0^{g,\phi}(\xi_n)\leq k$ as $m\rightarrow \infty$,
$y_0^m(\xi_n)\leq k$ holds for any $n$ and $m$.

For any fixed $m$, take $g_m=g +m\phi$, then $y_0^m(\xi_n) \rightarrow
y_0^m(\xi)$ as $n\rightarrow\infty$ and $y_0^m(\xi)\leq k$
holds for any $m$ according to  the continuous dependence
property of unconstrained BSDE.  For the fixed  $\xi\in L_T^2(\mathbb{R})$, the penalization method implies  $\mathcal{E}_0^{g,\phi}(\xi)\leq k$
and  $A_k$ is
closed under norm in $L_T^2(\mathbb{R})$. \hspace*{\fill}$\Box$

\begin{remark}\rm
In fact, Theorem 3.2 tells us the lower semi-continuity of $\mathcal{E}_0^{g,\phi}(\cdot)$. Since in some suitable conditions,  lower-semi-continuity is equivalent with the famous Fatou Property, thus by Remark 3.1, Theorem 3.2 can be thought as a corollary of Theorem 3.1 although we give an independent  proof.
\end{remark}

After discussions about the continuity  of $g_\Gamma$-solutions of CBSDEs in semi sense, we now come to study the full continuous properties.
We only consider the convex case, which means that $g$ and
$\phi$ are both convex.

\begin{pro}\label{con2}
 Suppose that $\phi(\omega,t,y,z)$ : $\Omega\times[0,T]\times \mathbb{R}\times \mathbb{R}^d\rightarrow \mathbb{R}^+$
and $g(\omega,t,y,z)$ : $\Omega\times[0,T]\times \mathbb{R}\times \mathbb{R}^d\rightarrow \mathbb{R}$ are
convex functions satisfying the assumptions $(A1)$ and $(A3)$, then  $D$ is convex in
$L_T^2(\mathbb{R})$.
\end{pro}

The proof of Proposition 3.1 is very similar to that of Propositions 2.1 and 2.2. So we omit it.

When both $g$ and $\phi$ are convex, Propositions 2.1 and 3.1 tell us that
$\mathcal{E}_0^{g,\phi}(\cdot)$ is  a  convex function  on
 $D$. This makes it possible for us to apply   wonderful results   in convex
 analysis to study the continuous dependence properties of CBSDEs with convex coefficients.
With the help of Theorem 4.1 stated in the appendix of this paper, we have

\begin{thm}
Suppose that  $\phi(\omega,t,y,z)$ : $\Omega\times[0,T]\times \mathbb{R}\times \mathbb{R}^d\rightarrow \mathbb{R}^+$
and $g(\omega,t,y,z)$ : $\Omega\times[0,T]\times \mathbb{R}\times \mathbb{R}^d\rightarrow \mathbb{R}$ are
convex functions satisfying the assumptions $(A1)$ and $(A3)$. If $A(4)$ holds, then
  $$E|\mathcal{E}_t^{g,\phi}(\xi_n)-\mathcal{E}_t^{g,\phi}(\xi)|^2\rightarrow
 0$$ when $E|\xi_n-\xi|^2\rightarrow 0$ as $n\rightarrow \infty$ for $\xi_n,\xi\in \mathring{D}$.
 \end{thm}
\proof Let $\overline{\xi}_n:=\xi_n\bigvee\xi$,
$\underline{\xi}_n:=\xi_n\bigwedge\xi$ for any $n$, then
$\{\overline{\xi}_n\}_{n=1}^\infty, \{\underline{\xi}_n\}_{n=1}^\infty$ are two sequences  of variables in $ L_T^2(\mathbb{R})$ and
converge  to $\xi$ in $L_T^2(\mathbb{R})$.

Let $\varphi(\eta):=E|\mathcal{E}_t^{g,\phi}(\eta)-\mathcal{E}_t^{g,\phi}(\xi)|^2 $ as in Theorem 3.1.  We want to  show $\alpha_n\rightarrow 0$ as $n\rightarrow\infty$
with $\alpha_n:=
E|\mathcal{E}_t^{g,\phi}(\overline{\xi}_n)-\mathcal{E}_t^{g,\phi}(\xi)|^2$.

 Let $K=\{\eta\in D|\eta\geq \xi \quad a.s.\}$  and $B_k=\{\eta\in K:\varphi(\eta)\leq k\}$. We first prove  $B_k$ is closed in $K$ and  obtain  the lower-semi continuity of $\varphi(\cdot)$ on $K$.
Suppose $\eta_n\in K, n=1,2,\cdots $ and $\eta_n\rightarrow\eta$ in norm, then $\varphi(\eta_n)=
E|\mathcal{E}_t^{g,\phi}(\eta_n)-\mathcal{E}_t^{g,\phi}(\xi)|^2\leq k$,  we need to show $\varphi(\eta)=
E|\mathcal{E}_t^{g,\phi}(\eta)-\mathcal{E}_t^{g,\phi}(\xi)|^2\leq k$. Suppose on the contrary $E|\mathcal{E}_t^{g,\phi}(\eta)-\mathcal{E}_t^{g,\phi}(\xi)|^2-k=\delta>0$. For  $\zeta=\eta,\xi,\eta_n$, we take
$\{y_t^m(\zeta),m=1,2,\cdots\}$ as in Proposition \ref{con1}.   Since $\varphi(\eta_n)=
E|\mathcal{E}_t^{g,\phi}(\eta_n)-\mathcal{E}_t^{g,\phi}(\xi)|^2\leq k$, we have
$$E|\mathcal{E}_t^{g,\phi}(\eta_n)-y^m_t(\xi)|^2\leq c$$
holds for any $n,m$ for some $c>0$.

Noting that
$y_t^m(\zeta)$ converges increasingly to
$\mathcal{E}_t^{g,\phi}(\zeta)$ as $m\rightarrow \infty$, we can find $M_1>0$ such that
$$E|\mathcal{E}_t^{g,\phi}(\xi)-y^m_t(\xi)|^2\leq \frac{\lambda^2\delta^2}{3(k+c)}\eqno(A)$$
for $m\geq M_1,  0<\lambda<1$.
At the same time, by dominated convergence theorem, we can find $M_2>0$ such that

$$E[y^m_t(\eta)-y^m_t(\xi)]^2-k>\lambda\delta\eqno(B)$$
when $m>M_2$.

Taking $M=\max\{M_1,M_2\}$, fix some $m>M$, we can find some $N(m)>0$ depending on $m$ such that

$$E[y^m_t(\eta_n)-y^m_t(\xi)]^2-k>\lambda\delta\eqno(C)$$
when $n>N$. By the  penalization method of CBSDEs and comparison property of BSDE, we have $\mathcal{E}_t^{g,\phi}(\eta_n)\geq y^m_t(\eta_n)\geq y^m_t(\xi)$, then

$$E[\mathcal{E}_t^{g,\phi}(\eta_n)-y^m_t(\xi)]^2-E|\mathcal{E}_t^{g,\phi}(\eta_n)-\mathcal{E}_t^{g,\phi}(\xi)|^2\geq E[\mathcal{E}_t^{g,\phi}(\eta_n)-y^m_t(\xi)]^2-k\geq E[y^m_t(\eta_n)-y^m_t(\xi)]^2-k  > \lambda\delta\eqno(D)$$
since $E|\mathcal{E}_t^{g,\phi}(\eta_n)-\mathcal{E}_t^{g,\phi}(\xi)|^2\leq k$.
On the contrary, by $H\ddot{o}lder's$ inequality,
$$E[\mathcal{E}_t^{g,\phi}(\eta_n)-y^m_t(\xi)]^2-E|\mathcal{E}_t^{g,\phi}(\eta_n)-\mathcal{E}_t^{g,\phi}(\xi)|^2\leq
\sqrt{3(k+c)}\sqrt{E|\mathcal{E}_t^{g,\phi}(\xi)-y^m_t(\xi)|^2}\leq \lambda \delta\eqno(E)$$
Thus we come to a contradiction of Equations (D) and (E) and the lower semi-continuity of the function $\varphi(\cdot)$ on $K$ is proved. With the convex
assumption of $g, \phi$, by Proposition \ref{con1} and  Proposition \ref{con2}, $\varphi(\cdot)$ is also convex on $K$.  With the help of some results in
convex analysis gathered in Appendix, we have $\alpha_n\rightarrow 0$ as $n\rightarrow \infty$.

Secondly, let $\beta_n:=
E|\mathcal{E}_t^{g,\phi}(\underline{\xi}_n)-\mathcal{E}_t^{g,\phi}(\xi)|^2.$
By Theorem \ref{cdp}, the continuous property from below,   we know that $\beta_n\rightarrow 0$ as
$n\rightarrow\infty$.

Finally, by Proposition 2.2, the dominated inequalities
$\mathcal{E}_t^{g,\phi}(\underline{\xi}_n)\leq\mathcal{E}_t^{g,\phi}({\xi}_n)\leq\mathcal{E}_t^{g,\phi}(\overline{\xi}_n)$
and
$$E|\mathcal{E}_t^{g,\phi}({\xi}_n)-\mathcal{E}_t^{g,\phi}(\xi)|^2\leq
\max\{E|\mathcal{E}_t^{g,\phi}(\overline{\xi}_n)-\mathcal{E}_t^{g,\phi}(\xi)|^2,
E|\mathcal{E}_t^{g,\phi}(\underline{\xi}_n)-\mathcal{E}_t^{g,\phi}(\xi)|^2\},$$
 guarantee our final result
 $$E|\mathcal{E}_t^{g,\phi}(\xi_n)-\mathcal{E}_t^{g,\phi}(\xi)|^2\rightarrow
 0\quad as \quad n\rightarrow \infty.$$ \hspace*{\fill}$\Box$

\begin{remark}\rm
  In El Karoui et al. \cite{kkppq}, the authors studied a new kind of BSDE reflected from below
by a barrier $S_t$ $(0\leq t\leq T)$. It is interesting that its solution  is just the
 minimal solution of BSDE with the constraint
$\phi(t,y_t,z_t)=(y_t-S_t)^-=0$ that was pointed out in Peng and Xu \cite{psgxmy2}. In this special case that $\{\xi\in L_T^2(\mathbb{R})|\xi\geq S_T, a.s.\}$ is  convex
 in $L_T^2(\mathbb{R})$, the existence of solution of RBSDE can be obtained by an a priori estimate and thus the minimal solution
is uniformly continuous with respect to the terminal value even when the coefficients are not convex.
\end{remark}

\begin{remark}\rm
 Under the
assumptions that $\mathcal{E}_0^{g,\phi}(\xi)<\infty$ for  $\xi\in
L_T^2(\mathbb{R})$ and  the constraint function $\phi(t,y,z)$ is the
distance function of $z/\sigma y$ from a convex closed set in $\mathbb{R}^n$,
  Karatzas and Shreve \cite{kiss},  and Civitanic et al. \cite{cks}  gave a representation of $\mathcal{E}_0^{g,\phi}(\xi)$ as a supremum  of a family  of  linear function
on $L_T^2(\mathbb{R})$, and  then, roughly speaking,
$\mathcal{E}_0^{g,\phi}(\xi)$    maybe uniformly continuous with
$\xi$ in some suitable normed space, even  the constraint function $\phi(t,y,z)$ is not convex.
 \end{remark}
 \begin{remark}\rm
By the comparison theorem of BSDEs, the increasing part $C_t$ of $g_\Gamma$-solution act  as some kind of  control to put up the solutions of BSDEs unidirectionally and thus make  the continuous properties hold only in semi sense generally as showed in the following example.
 \begin{example}\rm
 Let $g\equiv 0$ and $\Gamma_t=A_t\cup B_t, \,0\leq t\leq 1$, where
 $$A_t:=\{y:  0\leq y\leq 1\},\,0\leq t\leq 1$$
 and
  $$B_t:=\{y:  2-t\leq y\leq 2\},\,0\leq t\leq 1.$$
Obviously,  $\Gamma_t$ is  not convex.
  In this case, for any real number $\xi\in(1,2]$,
 $\mathcal{E}_0^{g,\phi}(\xi)=2$, but for $\xi=1$,  $\mathcal{E}_0^{g,\phi}(\xi)=1$.
It is obvious that $\mathcal{E}_0^{g,\phi}(\cdot)$-solution is not continuous at point $\xi=1$.
\end{example}
 \end{remark}

\section {Appendix: Some result about convex function}
Let $\varphi(\cdot)$ be a convex function defined on
$L_T^2(\mathbb{R})$, we assume  the effective domain of definition $D:=\{\xi\in
L_T^2(\mathbb{R})|-\infty<\varphi(\xi)<\infty\}$ of $\varphi(\cdot)$ is
closed in $L_T^2(\mathbb{R})$ and the interior of $D$, denoted by $\mathring{D}$,
is not empty. We omit the proofs of the following results since they can be founded in any book about convex analysis.
\begin{lem}
Suppose that $\varphi(\cdot):D\rightarrow \mathbb{R}$ is a  convex function and
$\xi_0$ is a point in $\mathring{D}$. If there exist a real number
$k>0$ and a neighborhood $O(\xi_0)\subset\mathring{D}$ of $\xi_0$
such that $\varphi(\xi)\leq k$ holds whenever  $\xi\in O(\xi_0)$,
then $\varphi(\cdot)$ is continuous at $\xi_0$.
\end{lem}

The next lemma tells us that if a convex function is bounded
in a neighborhood of some point $\xi_0\in \mathring{D}$, then it is
also locally bounded  at any $\xi\in\mathring{D}$.
\begin{lem}
Suppose  $\varphi(\cdot):D\rightarrow \mathbb{R}$ is a convex function
and bounded  in a neighborhood $O(\xi_0)\subset\mathring{D}$
of some $\xi_0\in\mathring{D}$, then for any $\xi\in
\mathring{D}$, there exist a neighborhood $O(\xi)\subset\mathring{D}$
of $\xi$ and a real number $k$ which dependents on $\xi$  such that
$\varphi(\eta)\leq k$ whenever $\eta\in O(\xi)$.
\end{lem}

 The following result is  well-known  in convex analysis which
ensures that the continuous property can be obtained by the lower semi-continuous property.
\begin{thm}
Suppose $\varphi(\cdot):D\rightarrow \mathbb{R}$ is a lower semi-continuous
convex function, then
 there exist some point $\xi_0\in\mathring{D}$ and a
neighborhood $O(\xi_0)\subset\mathring{D}$  such that
$\mathcal{E}_0^{g,\phi}(\xi)\leq k_0$ whenever $\xi\in O(\xi_0)$ for
some $k_0>0$, thus  it is continuous in $\mathring{D}$.
\end{thm}

\section*{Acknowledgements}
 The main work was done when the first author was a PhD student in Shandong University. He show
 his great thanks to Prof. Shige Peng for his valuable suggestions and providing stimulating working environment.


\begin{thebibliography}{}
 \bibitem{cks}Civitanic J., Karatzas I., Soner H.M.: Backward
stochastic differential equation with constraints on the
gains-process. Ann. Probab. 26(4), 1522--1551 (1998).


\bibitem{dhk}Drapeau S., Heyne G., Kupper M.: Minimal supersolutions of convex BSDEs. Ann. Probab. 41(6), 3697--4427   (2013).

\bibitem{kkppq}El Karoui N., Kapoudjian C., Pardoux E., Peng S.G. and Quenez M.C.: Reflected solution of backward SDE and related obstacle prolem for PDEs. Ann. Probab. 25(2), 702--737 (1997).

\bibitem{kpq}El Karoui N., Peng S.G., Quenez M.C.: Backward stochastic differential equations in Finance. Math. Finance. 7(1), 1--71 (1997).

\bibitem{fwz}Fan S.J., Wu Z.W., Zhu K.Y.: Continuous dependence properties on
solutions of backward differential equation. J. Appl. Math.
Comput.  24, 427--435 (2007).

\bibitem{fhsa}F\"{o}llmer H.,  Schied A.: Stochastic Finance. An introduction in discrete time. De Gruyter, Berlin. (2004).

\bibitem{jslpsg}Ji S.L., Peng S.G.: Terminal perturbation method for the backward approach to continuous time mean-variance portfolio selection. Stochastic Process. Appl. 118, 952--967 (2008).

\bibitem{kiss}Karatzas I.,  Shreve S.E.: Methods of mathematical finance. Springer Verlag New York. (1998).


\bibitem{pepsg} Pardoux E., Peng S.G.: Adapted solution of a backward stochastic differential equation. Systems \& Control Letters, 14, 55--62 (1990).

\bibitem{psg}Peng S.G.: Monotonic limit theorem of BSDE and nonlinear decomposition theorem of Doob-Meyer's type. Probab. Theory Relat. Fields. 113, 473--499 (1999).

\bibitem{psgxmy2}Peng S.G.,  Xu M.Y.: The smallest g-supermartingale and reflected BSDE
with single and double $L^2$ obstacles. Ann. I. H. Poincar\'{e}-PR
41, 605--630  (2005).

\bibitem{psgxmy1}Peng S.G., Xu M.Y.: Reflected BSDE with a constraint and its
applications in an imcomplete market. Bernoulli. 16(3),
614--640  (2010).



\bibitem{ger} Rosazza Gianin E.: Risk measures via
g-expectations. Insurance: Math. Economics 39, 19--34 (2006).

 \bibitem{yk}Yosida K.: Functional Analysis. Springer (1980).

\end{thebibliography}
\end{document}